\documentclass[11pt,twoside]{article}

\usepackage{amsfonts}
\usepackage{amssymb}
\usepackage{amsmath}
\usepackage{graphicx}
\usepackage{hyperref}
\numberwithin{equation}{section}
\newtheorem{Prop}{\bf Proposition}[section]
\newtheorem{Cor}{\bf Corollary}[section]
\newtheorem{defn}{\bf Definition}[section]
\newtheorem{Rem}{\bf Remark}[section]
\newtheorem{Ex}{\bf Example}[section]

\begin{document}
\def \b{\Box}

\begin{center}
{\Large {\bf On  almost groupoid associated to symmetric group $ S_{n} $}}

\end{center}

\begin{center}
{\large{\bf Mihai IVAN}}
\end{center}

\setcounter{page}{1}
\pagestyle{myheadings}

{\small {\bf Abstract}: The main purpose of this paper is to construct an example of a finite almost groupoid and namely: the almost groupoid associated to symmetric group $~S_{n}.$ 
Also, the remarkable properties of this are presented.}
{\footnote{{\it AMS classification:} 20L05, 20M10 18B40.\\
{\it Key words:} groupoid, almost groupoid, almost subgroupoid.}}

\section{Introduction}

The notion of groupoid from an algebraic point of view was first introduced  by H. Brandt in a 1926 paper \cite{brandt26}. Other equivalent definition
of groupoid and its properties appear in \cite{avmp20, ivan02, pata18}).

The Brandt groupoids were generalized by C. Ehresmann in \cite{ehre50}. After the introduction of topological and differentiable  groupoids by Ehresmann  in the 1950's, they
have been studied by many mathematicians with different approaches (\cite{brown06, ivan99}).
The topological groupoids and  Lie groupoids  play an essential role by their applications in analysis, differential geometry  and mechanics (\cite{guric14, giv5b, puta5a, miv03, miva05, miva06, mack87, rare01}) and so on.

Another approach to the notion of a groupoid is that of a structured groupoid.  This  concept is obtained by adding another algebraic structure
such that the composition of the groupoid is compatible with the operation of the added algebraic structure. The most important types of structured groupoids
 are  the following: group-groupoid  (\cite{brspen76}), vector space-groupoid (\cite{miva13}).\\
Furthermore,  topological groupoids, Lie groupoids and Lie algebroids are involved through their applications in the theory of classical and fractional dynamical systems 
(\cite{ igim11, giom13, opris5a, igmp11, ivmi16, mlag,  miva23a, migi18,  opris5b}.
 More details on basic properties and applications of groupoids, topological groupoids,
 Lie groupoids  and structured groupoids can be found  in the paper \cite{ivan24}.
 
A special field in groupoid theory is dedicated to the study of finite groupoids. In this context,   various  results are presented in the papers (\cite{beglpt21, miva02, miva03, miva26}).

In the paper \cite{miva23} was defined  the concept of almost groupoid. It is very important to note that the class of almost groupoids is included in the class of groupoids.

The paper is organized as follows. Section 2 is devoted to giving basic definitions and some results related to almost  groupoids. In Section 3 we give 
one interesting example of finite almost groupoids.\\[-0.7cm]

\section{Preliminaries on almost groupoids and its substructures}

Firstly, we will give the definition of an almost groupoid from a purely algebraic point of view.
\begin{defn}
(\cite{miva23}) {\rm   An {\it almost groupoid $ G$ over} $ G_{0}$  ({\bf in the sense of Brandt)} is a pair
$( G, G_{0} )$ of nonempty sets such that $~G _0\subseteq G,~ $ endowed with a surjective map
$\theta: G \to G_{0}$, a partially binary operation $~ m :G_{(2)}\rightarrow G ,~(x,y)\longmapsto  m ( x,y):=x\cdot y,~$ where\\
$~~~~~~~~~~~~~~~~~~~~~~~~~~~~ G_{(2)}:= \{ (x,y)\in G\times G | \theta(x) = \theta(y) \}$ \\
and a map $\iota:G \rightarrow G ,~~x\longmapsto \iota(x):=x^{-1},$ satisfying the following properties:\\
$({\bf AG1})~$ ({\it associativity})$~~(x\cdot y)\cdot z=x\cdot (y\cdot z)~$ in the sense that if
one of two products $ (x\cdot y)\cdot z $ and $ x\cdot (y\cdot z) $ is defined,
then the other product is also defined and they are equals;\\
$({\bf AG2})~$ ({\it units}): for each
$ x\in G~\Rightarrow~(\theta(x),x),~(x,\theta( x))\in G _{(2)}$
and  we have\\ $ \theta(x)\cdot x= x \cdot \theta(x)= x;$\\
$({\bf AG3})~$ ({\it inverses}): for each $ x\in G~\Rightarrow~(x, x^{-1}),~( x^{-1}, x) \in G_{(2)} $
 and we have\\ $ x\cdot x^{-1}= x^{-1}\cdot x= \theta(x).$}
\end{defn}

\markboth{M. Ivan}{On  almost groupoid associated to symmetric group $ S_{n}$}

If $ G $ is an almost groupoid over $ G_{0}, $  we will sometimes write $x\cdot y $ or $ x y$ for $ m(x,y). $ Also, the set $ G_{(2)} $ is called  the {\it set of composable pairs} of $ G. $

 An almost groupoid $ G $ over $ G_{0} $  with the \textit{structure functions} $ \theta
$ ({\it units map}),$ m $ ({\it multiplication}) and $ \iota $ ({\it inversion}), is denoted by $ (G, \theta,  m, \iota, G_{0}).~$  Whenever we write a product in a given almost groupoid, we are
assuming that it is defined.

In view of Definition 2.1, if $ x,y, z\in G, $ then:

 $(i)~~$ the product $~~x\cdot y\cdot z := (x\cdot y)\cdot z ~~$ {\it is defined} $~~~ \Leftrightarrow ~~~ \theta(x)=\theta (y) =\theta (z). $

$(ii)~~ \theta(x)\in G_{0} $ is {\it the unit} of $ x $ and  $ x^{-1}\in G $ is the {\it inverse} of $ x.$\\

For an almost groupoid $ G $  over $ G_{0}, $  we will use the following notations:  $ (G, \theta, m, G_{0} ) $
or $ (G, G_{0} ). $  Furthermore,  $ G_{0} $ is called the {\it units set} of $ G. $

The basic properties of the almost groupoids are given in the following proposition.
\begin{Prop} 
(\cite{miva23}) If $ (G,\theta,  m, \iota, G_{0}) $ is an almost groupoid, then:

$(i)~~~\theta(u)=u,~(\forall)u\in G_{0}.$

$(ii)~~u\cdot u=u ~$ and $~\iota (u)=u,~(\forall)u\in G_{0}.$

$(iii)~\theta(x\cdot y) =\theta (x), ~(\forall ) ( x,y ) \in G _{(2)}.$

$(iv)~\theta ( x^{-1}) =\theta ( x),~ (\forall) x\in G.$

$(v) $  If $ (x,y), (y,z)\in G_{(2)}~~\Rightarrow~~(x\cdot y, z), (x, y\cdot z) \in G_{(2)}~ $  and
$~(x \cdot y)\cdot z = x \cdot (y\cdot z).$

$(vi)~~ $  if $~(x,y)\in G_{(2)},$ then $~(y^{-1},x^{-1})\in
G_{(2)}~ $ and $~(x\cdot y)^{-1}=y^{-1}\cdot x^{-1}.$

{\it(vii)}$~~ (x^{-1})^{-1}=x, ~~(\forall) x\in G.$

{\it(viii)} $~~ x^{-1}\cdot(x\cdot y)=y~~$ and $~~(x\cdot y)\cdot y^{-1}=x,$ for all $~(x,y)\in G_{(2)}.$
\end{Prop}
\begin{Cor} If $ (G,\theta,  m, \iota, G_{0}) $ is an almost groupoid, then:

$(i)~~~~~~ \theta(\theta(x))= \theta(x), ~~(\forall) x\in G. $

$(ii)~~~~$  If $~ x,y\in G,~$  then:  $~~ x\cdot y ~~$ is defined $~~~\Leftrightarrow~~~ y\cdot x~~$  is defined.

$(iii)~~~ (\forall)~ a\in G~~\Rightarrow ~~$ the products $~ a^{2}:= a\cdot a, ~ a^{3}:= a^{2}\cdot a~ $ are defined.
\end{Cor}
{\bf Proof.} $(i)~$ For each $ x\in G, $ we have $~ \theta(\theta(x)) \stackrel{(AG3)}{=}~\theta (x\cdot x^{-1}) \stackrel{(2.1(iii))}{=}~  \theta(x).$

$(ii)- (iii)~$ These statements follow from definition and Proposition 2.1.\hfill$\Box$
\begin{Prop} Let $ (G,\theta,  m, \iota, G_{0}) $  be an almost groupoid. For each $ u \in G_{0}, $ the set $ G(u):=\theta^{-1}(u)= \{ x\in G | \theta(x)= u \} $ is
a group.
\end{Prop}

The group $ G(u) $ for $ u\in G_{0} $ is called the {\it isotropy group} of $ G $ at $ u.$

An {\it almost groupoid $ G$ over} $ G_{0}$ is called {\it abelian} or {\it commutative} if the isotropy group $ G(u) $ is abelian for all $ u\in G_{0}. $
\begin{Ex}
{\rm $(i)~$ A group $ G $ having $ e $ as unity, is an almost groupoid over   $ \{e
\}$ with respect to structure functions: $~\theta (x): = e, ,(\forall)x\in G; ~G_{(2)}= G\times G,~ m (x,y):= xy,~ (\forall) x,y\in G
  $ and $ \iota:G \to G,~\iota(x):= x^{-1}, (\forall)x\in G.$

$(ii)~$ A nonempty set $ G_{0} $ may be regarded to be an
almost groupoid over $ G_{0}, $ called the {\it null almost groupoid} associated
to $ G_{0}. $ For this, we take $~ \theta=\iota = Id_{G_{0}} $ and\\
 $ u\cdot u = u, ~(\forall)  u\in G_{0}.$}
\end{Ex}

Next,  we present some substructures in almost groupoids and recall  several properties of them which generalize well-known results in group theory.

\begin{defn}
{\rm  Let $(G ,\theta, m, \iota, G_{0})$ be an almost groupoid. A pair of nonempty subsets $(H,H_{0})$
where $ H\subseteq G $ and ${H_0}\subseteq G_{0} $, is called
\textit{almost subgroupoid} of $ G, $ if:

$(1.1)~~\theta (H)=H_{0};$

$(1.2)~~H $ is closed under multiplication and inversion,
that is:\\[0.1cm]
$(1.2.1)~(\forall)~x,y\in H$  such that $(x,y)\in G_{(2)}
\Longrightarrow ~x\cdot y\in H;$\\[0.1cm]
 $~(1.2.2)~(\forall)~x\in
H\Longrightarrow x^{-1}\in H.$}
\end{defn}
\begin{defn}
{\rm Let $(H,H_{0})$ be an almost subgroupoid of the almost groupoid $ (G,G_{0}).$

{\it(i)} $ (H, H_{0}) $ is said to be a {\it wide almost subgroupoid} of $ (G, G_{0} ), $ if $~H_{0}= G_{0},$   that is  $ H $  and $ G $  have the same units.

{\it(ii)} A wide almost subgrupoid $ (H,H_{0})$  of $ (G,G_{0}),$  is called {\it normal almost subgroupoid,} if for all $ g\in G $ and for all
 $ h \in H $ such that the product $ g\cdot h\cdot {g^ {-1}} $ is defined, we have $ g\cdot h\cdot {g^ {-1}}\in H .$}
\end{defn}
\begin{Ex}
{\rm  $~(i)~$ If $~(G, G_{0})~$ is an almost groupoid, then $ G_{0} $ is a wide almost subgroupoid of $ G. $ \\
$~(ii)~$ $ G_{0} $ is a normal almost subgroupoid of the almost groupoid  $ (G, G_{0}). $  
}\end{Ex}

\begin{Prop} (\cite{miva23}) Let $ (G, G_{0}) $  be an almost groupoid.  Then:\\
$(i)~~$ The isotropy group  $~G(u), u\in G_{0}~$ is an almost subgroupoid of $ G. $\\
$(ii)~~$  The set $~Is(G) = \cup_{u\in G_{0}}G(u)\subset G~$ is a normal almost subgroupoid of $~G.~$
\end{Prop}

The almost groupoid $~(Is(G), G_{0}) $ is called the {\it isotropy almost subgroupoid} of $~G. $

Let  $~(G, \theta, m, \iota, G_{0}) $  be an almost groupoid and $ a\in G.$   We define the set $ <a> $, given by:\\[-0.4cm]
\begin{equation}
\begin{array}{l}
< a> =\{ a^{n} ~|~(\forall)~n\in {\bf Z}\},~\hbox{where}~~ a^{0}=\theta (a),~~a^{n}= a^{n-1}\cdot a,~~a^{-n}=(a^{-1})^{n},~n\in {\bf N}^{\ast}.\\
\end{array}\label{2.1}
\end{equation}
It is easy to prove that the set $<a>$ is an almost sugbroupoid,  called the {\it cyclic almost subgroupoid generated by} $ a\in G. $
\begin{Rem}~{\rm Let $~(G, G_{0}~$ be an almost groupoid  and $~a\in G.$  Then:\\
$(i)~~~ < a >\subset G(\theta(a))~~$ and $~<a>~$ is a subgroup of the isotropy group $~G(\theta(a)). $\\
$(ii)~$ If $~(\exists)~ n\in {\bf N}, ~n\geq 2~ $ such that $~a^{n} =\theta(a),$ then  $ < a > $ is a cyclic subgroup of order $ n~ $ and
$~<a>  =\{ \theta (a), a, a^{2}, ..., a^{n-1} \}.$}
\end{Rem}
\begin{Rem}~{\rm For definitions and basic results related to almost groupoid morphisms, see \cite{ ivan25, miva23}.}\\[-0.5cm]
\end{Rem}

\section{The almost groupoid associated to symmetric group $~S_{n} $}

Symmetric groups are some of the most essential types of finite groups. The symmetric group $ S_{n} $ defined over a finite set of $ n $ elements consists of the permutations that can be performed on the elements. The
order (number of permutations) of the symmetric group $ S_{n} $ is $~n!. $ Usually the elements of $ S_{n} $ are given by bijective function $ \sigma: \{ 1, 2, . . ., n\} \rightarrow \{ 1, 2, . . ., n\}. $  The group operation on $ S_{n} $ is composition functions. It is 
well-known that $ S_{n} $ is non commutative for $ n\geq 3. $ Denote the elements of $ S_{n} $ by $\sigma_{i}, ~ i=\overline{ 1, n!},$ where $\sigma_{1} = Id $ (identity permutation).

The set of all even permutations in $ S_{n} $ is denoted by $ A_{n} $ and is called the {\it alternating group}. $~A_{n}$ is a normal subgroup of $~S_{n}.$  It is a finite group of order $ \frac{n!}{2}.$

Let $ M=\{ 0, 1\}~$ a set consisting of two elements. We now define an almost groupoid  structure on the set $ G:=M\times S_{n}= \{ (a, \sigma_{i}) ~|~ a\in M,~ \sigma_{i}\in S_{n}~\hbox{for}~i=\overline{1,n!}\}.$

Let $ G_{0} = \{ (a,\sigma_{1}) | a\in M\}. $  Then $~ (G, \theta, \odot, \iota, G_{0})~$ is an almost groupoid over $ G_{0}, $ where  the set of
composable pairs $~ G_{(2)}~$ and  its structure functions are given by:\\
$\theta : G=M\times S_{n}~ \rightarrow~ G_{0}, (a,\sigma)~\rightarrowtail \theta (a,\sigma):=(a,\sigma_{1});$\\
$ G_{(2)} = \{ ((a,\sigma_{i}), (b, \sigma_{j}))\in G^{2}~|~\theta (a,\sigma_{i}) =\theta (b, \sigma_{j})\} = \{ ((a,\sigma_{i}), (b, \sigma_{j}))\in G^{2}~|~ a=b \};$\\
$(a,\sigma_{i})~\odot~(a,\sigma_{j}):=(a, \sigma_{i}\circ \sigma_{j}),~ (\forall)~a \in M~, \sigma_{i}, \sigma_{j}\in S_{n}~$ and\\
$\iota(a,\sigma_{i}):=(a, \sigma_{i}^{-1}), (\forall)~a\in M,~\sigma_{i}\in S_{n}. $

It is easy to verify that the conditions $ ({\bf AG1}) - {\bf AG3}) $ of Definition 3.1 hold. First, we consider $~ x, y, z\in G ,~$ where $~x=(a, \sigma_{i}), y=(b, \sigma_{j}), z=(c, \sigma_{k}). $  The product $ x\odot y \odot z $ is defined if and only if $~a = b = c.$ We have:\\
$ (1)~~~(x\odot y)\odot z= ((a, \sigma_{i})\odot (a, \sigma_{j}))\odot (a, \sigma_{k})= (a, \sigma_{i}\circ \sigma_{j})\odot (a, \sigma_{k})=\\
 (a, \sigma_{1}\circ \sigma_{j}\circ \sigma_{k}) =x\odot (y\odot z); $\\
$ (2)~~~x\odot \theta(x)= (a, \sigma_{i})\odot \theta(a, \sigma_{i}))= (a, \sigma_{i})\circ (a, \sigma_{1})=(a, \sigma_{i}\circ \sigma_{1})= (a, \sigma_{i})= x =\theta(x)\odot x;$\\
$ (3)~~~x\odot \iota(x)= (a, \sigma_{i})\odot \iota(a, \sigma_{i}))= (a, \sigma_{i})\odot (a, \sigma_{i}^{-1})= (a, \sigma_{i}\circ \sigma_{i}^{-1})= (a, \sigma_{1})=\\
 \theta(x)= \iota(x)\odot x.$

We conclude that $ G = M\times S_{n}~$ is an almost groupoid over $ G_{0}=\{ (0, \sigma_{1}), (1, \sigma_{1}).~$  Observe that $~ (G, G_{0})~$ is a finite almost groupoid of order $~2 n!.$ This almost groupoid is called the {\it almost groupoid associated to $ S_{n} $}.

Let be the set $ H:=M\times A_{n}\subset M\times S_{n}.~$ It is easy to prove that $ (M\times A_{n}, \theta, \odot, \iota,  G_{0}) $ is a normal almost subgroupoid of $ G. $ It is called the {\it almost groupoid associated to}  $ A_{n}. $ The order of the almost groupoid $ M\times A_{n} $ is $ n!. $

We now consider the symmetric group $~S_{3} =\{ \sigma_{i} | i=\overline{1,6} \},~$ where:\\
$~~~ \sigma_{1}=\left(\begin{array}{ccc}
1 & 2 & 3\\
1 & 2 & 3
\end{array}\right)=Id,~ ~ \sigma_{2}=\left (\begin{array}{ccc}
1 & 2 & 3\\
2 & 3 & 1
\end{array}\right),~ \sigma_{3}=\left (\begin{array}{ccc}
1 & 2 & 3\\
3 & 1 & 2
\end{array}\right),\\~ \sigma_{4}=\left (\begin{array}{ccc}
1 & 2 & 3\\
2 & 1 & 3
\end{array}\right),~ \sigma_{5}=\left (\begin{array}{ccc}
1 & 2 & 3\\
3 & 2 & 1
\end{array}\right),~ \sigma_{6}=\left (\begin{array}{ccc}
1 & 2 & 3\\
1 & 3 & 2
\end{array}\right).~$

$\bullet~~$ In the following we refer to the almost groupoid $ G=M\times S_{3} $ associated to the symmetric group $ S_{3}. $ In this case,$~G=\{ u_{1}, u_{2}, g_{1}, g_{2}, . . ., g_{9}, g_{10} \}, $   where:\\
$~u_{1}=(0,\sigma_{1}),~ u_{2}=(1,\sigma_{1}), $\\
$g_{1}=(0,\sigma_{2}),~ g_{2}=(0,\sigma_{3}),~g_{3}=(0,\sigma_{4}),~g_{4}=(0,\sigma_{5}),~g_{5}=(0,\sigma_{6}),$\\
$g_{6}=(1,\sigma_{2}),~g_{7}=(1,\sigma_{3}),~g_{8}=(1,\sigma_{4}),~g_{9}=(1,\sigma_{5}),~g_{10}=(1,\sigma_{6}). $

{\small
$$\begin{array}{|r|c|c|c|c|c|c|c|c|c|c|c|c|} \hline
       g    & u_{1} & u_{2} & g_{1} & g_{2} & g_{3} & g_{4} & g_{5} & g_{6} & g_{7} & g_{8} & g_{9} & g_{10} \\ \hline
  \theta(g) & u_{1} & u_{2} & u_{1} & u_{1} & u_{1} & u_{1} & u_{1} & u_{2} & u_{2} & u_{2} & u_{2} & u_{2} \cr \hline
     \iota(g)  & u_{1} & u_{2}  & g_{2} & g_{1} & g_{3} & g_{4} & g_{5} & g_{7} & g_{6} & g_{8} & g_{9} & g_{10}\cr \hline
\end{array}$$
$$\begin{array}{|r|c|c|} \hline
       u    & u_{1} & u_{2} \\ \hline
     G(u)  & \{u_{1}, g_{1},  g_{2}, g_{3}, g_{4}, g_{5}\} & \{ u_{2}, g_{6}, g_{7},  g_{8}, g_{9},g_{10} \}\cr \hline
\end{array}$$}

We now give some examples for completing the above tables:\\
$\theta(u_{2}) =\theta(1,\sigma_{1})=(1,\sigma_{1})=u_{2};~~~~~~~~~~~~~~~\theta(g_{5})=\theta(0,\sigma_{6})= (0, \sigma_{1})=u_{1}.$\\
$\iota(u_{1})=\iota(0,\sigma_{1})=(0, \sigma_{1}^{-1})=(0, \sigma_{1})= u_{1};~~
\iota(g_{5})=\iota(0,\sigma_{6})=(0, \sigma_{6}^{-1})=(0,\sigma_{6})=g_{5}.$

The multiplication operation $~\odot $ on  $~G= M\times S_{3}~$ is given in the following table:
{\small
$$\begin{array}{|r|c|c|c|c|c|c|c|c|c|c|c|c|} \hline
\odot  & u_{1} & u_{2} & g_{1} & g_{2} & g_{3} & g_{4} & g_{5} & g_{6} & g_{7} & g_{8} & g_{9} & g_{10}\\ \hline
u_{1}  & u_{1} &  & g_{1} & g_{2} & g_{3} & g_{4} & g_{5} &  &   &  &  &  \cr \hline
u_{2}  &  & u_{2}  &  &  &  &  &  & g_{6} & g_{7}  & g_{8} & g_{9} & g_{10} \cr \hline
g_{1}  & g_{1} &  & g_{2} & u_{1} & g_{4} & g_{5} & g_{3} &  &  &   &  &   \cr \hline
g_{2}  & g_{2} &  & u_{1} &  g_{1} & g_{5} & g_{3} & g_{4} &  &  &  &  & \cr \hline
g_{3}  & g_{3} &  & g_{5} & g_{4} & u_{1} & g_{2} & g_{1} &  &   &  &  &    \cr \hline
g_{4} & g_{4} & & g_{3} & g_{5} & g_{1} & u_{1} & g_{2} &  & & &  &  \cr \hline
g_{5}  & g_{5}  &  & g_{4} & g_{3} & g_{2} & g_{1} &  u_{1} &  & & &  &    \cr \hline
g_{6} & & g_{6} &  &  &  & & & g_{7} & u_{2} & g_{9} & g_{10} & g_{8}  \cr \hline
g_{7}  &  & g_{7} &  &  &  &  &  & u_{2} & g_{6} & g_{10} & g_{8} & g_{9}   \cr \hline
g_{8} &  & g_{8} &  &  & &  & & g_{10} & g_{9} & u_{2} & g_{7} & g_{6}   \cr \hline
g_{9}  &  & g_{9}  & &  & & &  & g_{8} & g_{10} & g_{6} & u_{2} & g_{7} \cr \hline
g_{10} &  & g_{10} &  & &  &  & &  g_{9} & g_{8} & g_{7} & g_{6} & u_{2}  \cr \hline
\end{array}.~$$}

For to illustrate the calculation of the elements $ x \odot y $ for $ x,y \in G ~$ we give the following examples:\\
$u_{1}\odot g_{4}= (0,\sigma_{1})\odot (0,\sigma_{5})= (0, \sigma_{1}\circ \sigma_{5})=(0,\sigma_{5})=g_{4};$\\
$g_{2}\odot  g_{4}= (0,\sigma_{3})\odot (0,\sigma_{5})= (0, \sigma_{3}\circ \sigma_{5})=(0,\sigma_{4})=g_{3};$\\
$g_{7}\oslash g_{7}= (1,\sigma_{3})\odot (1,\sigma_{3})= (1, \sigma_{3}\circ \sigma_{3})=(1,\sigma_{2})=g_{6}.$

The absence of the element from the ligne  $~l_{i}~$ and the column $~c_{j}~$ in the table of composition law $~\odot~$ indicates the fact that the pairs $~x=(a,\sigma_{i})\in G~$ and $~y=(b,\sigma_{j})\in G~ $ are not composables. For example, the element $~g_{2}\odot g_{8} ~$ is not defined, since $ \theta(g_{2})=u_{1},~ \theta(g_{8})= u_{2}~$  and $~u_{1}\neq u_{2}. $

It is easy to observe that  $ Is(G) = G .$  Also,  $ G = M\times S_{3} $ is
a non commutative almost groupoid.

The cyclic almost subgroupoid generated by $ g_{1}= (0,\sigma_{2} ) $ is $ <g_{1}> =\{ u_{1}, g_{1}, g_{2}\}, $  since $~g_{1}^{2} = g_{1}\odot g_{1}=g_{2} ~$ and  $~g_{1}^{3} = g_{2}\odot g_{1}=u_{1}.~$ Also,
the cyclic almost subgroupoid generated by $ g_{9}= (1,\sigma_{5} ) $ is $ <g_{9}> =\{ u_{2}, g_{9}\}, $  since $~g_{9}^{2} = g_{9}\odot g_{9}=u_{2}.~$

$\bullet~~$ We now refer to the almost groupoid $ H=M\times A_{3} $ associated to  group $ A_{3}. $ In this case,$~H=\{ u_{1}, u_{2}, g_{1}, g_{2}, g_{6}, g_{7} \} $
is a finite almost groupoid of order $ 6, $ where:
{\small
$$\begin{array}{|r|c|c|c|c|c|c|} \hline
       g    & u_{1} & u_{2} & g_{1} & g_{2}& g_{6} & g_{7}\\ \hline
  \theta(g) & u_{1} & u_{2} & u_{1} & u_{1} & u_{2} & u_{2} \cr \hline
     \iota(g)  & u_{1} & u_{2}  & g_{2} & g_{1} & g_{7} & g_{6}  \cr \hline
\end{array},~~~~~
\begin{array}{|r|c|c|} \hline
       u    & u_{1} & u_{2} \\ \hline
     G(u)  & \{u_{1}, g_{1},  g_{2}\} & \{ u_{2}, g_{6}, g_{7} \}\cr \hline
\end{array}$$}
{\small
$$\begin{array}{|r|c|c|c|c|c|c|} \hline
\odot  & u_{1} & u_{2} & g_{1} & g_{2} & g_{6} & g_{7} \\ \hline
u_{1}  & u_{1} &  & g_{1} & g_{2} &  &  \cr \hline
u_{2}  &  & u_{2}  &  &   & g_{6} & g_{7}  \cr \hline
g_{1}  & g_{1} &  & g_{2} & u_{1} & &   \cr \hline
g_{2}  & g_{2} &  & u_{1} &  g_{1} &  & \cr \hline
g_{6} & & g_{6} &  &   & g_{7} & u_{2}  \cr \hline
g_{7}  &   & g_{7} &  & & u_{2} & g_{6}  \cr \hline
\end{array}.~$$}

Finally, the almost groupoid  $ H=M\times A_{3} $ is a normal almost subgroupoid of $ G=M\times S_{3}. $
\begin{Rem}
{\rm The almost groupoid $~M\times S_{3} $ can be interpreted as an assembly consisting of a bar on which two points (denoted by $0$ and $ 1$) are fixed so as to hang two equilateral triangles that can rotate independently or simultaneously.}
\end{Rem}

Author's adress\\[0.2cm]
\hspace*{0.7cm}West University of Timi\c soara. Seminarul de Geometrie \c si Topologie\\
\hspace*{0.7cm}  Teacher Training Department\\
\hspace*{0.7cm} Bd. V. P{\^a}rvan,no.4, 300223, Timi\c soara, Romania\\
\hspace*{0.7cm}E-mail: mihai.ivan@e-uvt.ro\\

\end{document}